\documentclass[english,11pt]{article}
\usepackage{lmodern}
\usepackage[T1]{fontenc}
\usepackage[utf8]{inputenc}
\usepackage[pdftex]{graphicx}
\usepackage{verbatim}
\usepackage[width=\textwidth]{caption}
\usepackage{amsmath,bbm,amssymb,amsopn,bm,mathtools,dsfont}
\usepackage{a4wide}
\usepackage[english]{babel}
\usepackage{float,placeins,flafter,array}
\usepackage{color}
\usepackage{url}
\usepackage{xr}
\usepackage{natbib}
\usepackage{marginnote}

\newtheorem{thm}{Theorem}

\newtheorem{lem}{Lemma}

\newtheorem{defn}{Definition}
\newtheorem{rem}{Remark}

\newcommand{\card}{\text{card}}
\newcommand{\ind}[1]{\mathds{1}\left\{#1\right\}}

\newcommand{\E}{\mathbb E}

\newcommand{\V}{\mathbb V}

\newcommand{\I}{\mathbb I}
\newcommand{\N}{\mathbb N}

\renewcommand{\j}{\bm{j}}

\newcommand{\un}{\boldsymbol{1}}

\newcommand{\Cov}{\mathbb{C}ov}
\newcommand{\e}{\bm{e}}

\newcommand{\zero}{\boldsymbol{0}}

\renewcommand{\i}{\bm{i}}

\newcolumntype{C}[1]{>{\centering\let\newline\\\arraybackslash\hspace{0pt}}m{#1}}


\date{}

\begin{document}

\title{The Marcinkiewicz-Zygmund law of large numbers for exchangeable arrays}
\author{Laurent Davezies\thanks{CREST-ENSAE, France. laurent.davezies@ensae.fr.} \and Xavier D'Haultf\oe uille\thanks{Corresponding author. CREST-ENSAE, 5 avenue Henry Le Chatelier, 91120 Palaiseau, France. xavier.dhaultfoeuille@ensae.fr.} \and Yannick Guyonvarch\thanks{INRAE, France. yannick.guyonvarch@gmail.com.}}
\maketitle

\begin{abstract}
 We show a Marcinkiewicz-Zygmund law of large numbers for jointly, dissociated exchangeable arrays, in $L^r$ ($r\in (0,2)$) and almost surely. Then, we obtain a law of iterated logarithm for such arrays under a weaker moment condition than the existing one. \\

\end{abstract}

\section{Introduction}

Random variables with multidimensional indexes are useful to model interactions between several units of a population. An example is international trade, where each observation corresponds to a pair of countries, one exporting and the other importing. Another example is network formation, where binary variables between pairs of units (e.g., individuals) measure whether they share a link or not. Multi-indexed random variables are also useful in case of multiway clustering. For instance, wages of individuals may be indexed by geographical areas and sectors of activity, reflecting potential dependence in wages within the same area or sector.

\medskip
There are two important properties for such arrays of random variables. The first, dissociation, was introduced by \cite{mcginley1975}. It generalizes mutual independence, by allowing for dependence when elements of the array share at least one common component in their respective indexes. The second, joint exchangeability, was introduced by \cite{silverman1976limit}, and is very similar to exchangeability for sequences. Taken together, these two conditions turn out to be sufficient to establish strong laws of large numbers, central limit theorems and weak convergence of empirical processes under the exact same moment conditions as in the i.i.d. case. We refer respectively to \cite{Eaglesonweber78}, \cite{silverman1976limit} and \cite{davezies2021empirical} for proofs of these statements.

\medskip
The aim of this note is to show that we can also extend the Marcinkiewicz-Zygmund law of large numbers to dissociated, jointly exchangeable arrays. Recall that this law extends the usual law of large numbers by relating the convergence rate of partial sums to moment conditions on the variables.  
As with i.i.d. data, our result holds both in $L^r$ ($r\in (0,2)$) and almost surely. The proof in $L^r$ for the more difficult case ($r>1$) is based on a symmetrization lemma established in \cite{davezies2021empirical}. Almost-sure convergence follows using the so-called Aldous-Hoover-Kallenberg (AHK) representation \citep{Aldous1981,Hoover1979,kallenberg1989} and a similar strategy as that used for U-statistics \citep{Gine1992}. As a consequence of our result, we obtain a law of iterated logarithm for dissociated, jointly exchangeable arrays with a finite second moment. We thus improve upon \cite{scott1985law}, who obtain the same result under the existence of a moment of order $2+\delta$ for some $\delta>0$.

\medskip
We first present in Section \ref{sec:res} the set-up and the main results. Section \ref{sec:proof} is devoted to the proof of the main theorem in the difficult case ($r>1$). The supplementary material includes in particular the proofs of the easier part ($0<r<1$) and of the law of iterated logarithm.

\section{The set-up and results}
\label{sec:res}

We first define formally dissociated and jointly exchangeable arrays. Let $\N^+=\N\backslash\{0\}$, $k\in \N^+$,  $\card(A)$ denote the cardinal of a set $A$ and for any $A\subset \N^k$,
$$\overline{A}=\left\{(a_1,...,a_k)\in A: \card(\{a_1,...,a_k\})=k\right\}.$$
Then, we denote by $\mathbb{I}_k=\overline{\N^{+k}}$ the set of $k$-tuples of $\mathbb{N}^+$ without repetition. Similarly, for any $n\in \N^+$, we let $\I_{n,k}=\overline{\{1,...,n\}^k}$. For any $A\subset \N^+$, we use $\mathfrak{S}(A)$ to denote the set of permutations on $A$. For any $\i=(i_1,...,i_k)\in \N^{+k}$ and $\pi\in \mathfrak{S}(\N^+)$, we let $\pi(\i)=(\pi(i_1),...,\pi(i_k))$. Finally, we let $\bm{k}=(1,...,k)$.

\begin{defn}[Dissociated and exchangeable arrays]
\label{def:dgp}~~\\
$X:=(X_{\i})_{\i \in \I_k}$ is a dissociated array if for any $A,B$ disjoint subsets of $\mathbb{N}^+$ with $\min(\card(A),$ $\card(B))\geq k$, $(X_{\i})_{\i \in \overline{A^k}}$ is independent of $(X_{\i})_{\i \in \overline{B^k}}$.\\
$X$ is a jointly exchangeable array if for any $\pi\in \mathfrak{S}(\mathbb{N}^+)$, $X\overset{d}{=}(X_{\pi(\i)})_{\i \in \I_k}$.
\end{defn}

For $k=1$, a jointly exchangeable array is an exchangeable sequence of random variables whereas a dissociated jointly exchangeable array is a sequence of i.i.d. random variables. Otherwise ($k>1$), jointly exchangeable arrays may exhibit more complex forms of dependence: with $k=2$, $X_{(i_1,i_2)}$ and $X_{(i_1,j_2)}$, but also $X_{(i_1,i_2)}$ and $X_{(j_2,i_1)}$ may be dependent, for instance. In the international trade example, exports from China to the USA may be correlated with exports from the USA to France, simply because the USA are a large country.

\medskip
Note that with multiway clustering, the different components of the index (e.g., areas and sectors in the wage example above) do not belong to the same population. Separate exchangeability, where the permutation used may differ for all components of the index, is then better suited to such data than joint exchangeability. However, joint exchangeability is more general so the result below also applies to separately exchangeable arrays.

\begin{thm}\label{thm:MZL}
Let $X$ be a jointly exchangeable array. Then:
\begin{enumerate}
	\item If $\E(|X_{\bm{k}}|^r)<\infty$ for $1\leq r<2$ and $X$ is dissociated, $$\frac{1}{n^{k-1+1/r}}\sum_{\i\in \I_{n,k} }\left(X_{\i}-\E(X_{\bm{k}})\right)\rightarrow 0 \text{ in } L^r \text{ and almost surely (a.s.)}.$$
	\item If $\E(|X_{\bm{k}}|^r)<\infty$ for $0<r<1$, $n^{-k/r}\sum_{\i \in \I_{n,k}}X_{\i}\rightarrow 0$  in $L^r$ and a.s.
\end{enumerate}
\end{thm}

The  case $0<r<1$ follows by simply noting that existing proofs \citep[see][]{gut2013probability,Gine1992} in fact only rely on the variables being identically distributed. The case $r=1$ is Theorem 3 in \cite{Eaglesonweber78}. The case $r>1$ is not as straightforward however. The key ingredients for convergence in $L^r$ are a symmetrization lemma established in \cite{davezies2021empirical}, see Lemma \ref{lem:sym} below, and the Khintchine inequality. Almost-sure convergence follows using a similar strategy as for U-statistics. To understand the similarities and differences with U-statistics,  let us consider the case $k=3$. By the AHK representation of dissociated and jointly exchangeable arrays, there exists a function $\tau$ and  i.i.d. and  uniform $[0,1]$ variables $(U_{\{i\}})_{i\geq 1}$, $(U_{\{i,i_2\}})_{(i,i_2)\in \I_{2}}$ and  $(U_{\{i,i_2,i_3\}})_{(i,i_2,i_3)\in \I_3}$ such that for all $\i=(i_1,i_2,i_3)\in\I_3$,
$$X_{\i}=\tau\left(U_{\{i_1\}},U_{\{i_2\}},U_{\{i_3\}},U_{\{i_1,i_2\}},U_{\{i_1,i_3\}},U_{\{i_2,i_3\}},U_{\{i_1,i_2,i_3\}}\right).$$
This result shows that $X_{\i}$ is close to, but distinct from, a U-statistic, where the last four variables would not appear. Accordingly, we use a decomposition akin to but different from the Hoeffding decomposition in the proof. Specifically, let us define
\begin{align*}
	H_0(X)_{\i}=&\E(X_{\i}),\\
	H_1(X)_{\i}=&\E\left(X_{\i}|U_{\{i_1\}}\right)+\E\left(X_{\i}|U_{\{i_2\}}\right)+\E\left(X_{\i}|U_{\{i_3\}}\right)-3\E(X_{\i}),\\
	H_2(X)_{\i} =&Q_{(1,1,0)}(X)_{\i}+Q_{(1,0,1)}(X)_{\i}+Q_{(0,1,1)}(X)_{\i},\\
	H_3(X)_{\i}=&X_{\i}-\E\left(X_{\i}|U_{\{i_1\}},U_{\{i_2\}},U_{\{i_1,i_2\}}\right)-\E\left(X_{\i}|U_{\{i_1\}},U_{\{i_3\}},U_{\{i_1,i_3\}}\right)\\
	&-\E\left(X_{\i}|U_{\{i_2\}},U_{\{i_3\}},U_{\{i_2,i_3\}}\right)+\E\left(X_{\i}|U_{\{i_1\}}\right)\\
	&+\E\left(X_{\i}|U_{\{i_2\}}\right)+\E\left(X_{\i}|U_{\{i_3\}}\right)-\E\left(X_{\i}\right),
\end{align*}
where $Q_{(1,1,0)}(X)_{\i}=\E\left(X_{\i}|U_{\{i_1\}},U_{\{i_2\}},U_{\{i_1,i_2\}}\right)-\E\left(X_{\i}|U_{\{i_1\}}\right)-\E\left(X_{\i}|U_{\{i_2\}}\right)+\E(X_{\i})$ and $Q_{(1,0,1)}(X)_{\i}$ and $ Q_{(0,1,1)}(X)_{\i}$ are defined similarly. Then,
\begin{equation}
X_{\i}= H_0(X)_{\i}+H_1(X)_{\i}+H_2(X)_{\i}+H_3(X)_{\i},
	\label{eq:decomp2}
\end{equation}
where the four terms can be shown to be orthogonal in $L^2$. Decomposition \eqref{eq:decomp2} allows us to extend \cite{Gine1992}'s approach to our set-up.

\begin{rem}
The proof shows that a result similar to Point 1 holds without dissociation. The only difference, then, is that one should replace $\E(X_{\bm{k}})$ by $\E(X_{\bm{k}}|U_\emptyset)$ in the statement of the theorem, where $U_\emptyset$ is an additional argument in the AHK decomposition.
\end{rem}

\begin{rem}
	For $k\geq 2$, the function $r\mapsto (k-1+1/r) \ind{r<1} + k/r \ind{r\geq1}$ is continuous but not differentiable at $r=1$. There is therefore a kink at $r=1$ in the rate of convergence. 
\end{rem}

\begin{rem}
The result is sharp in the following sense: there exist exchangeable arrays $(X_{\i})_{\i\in\I_k}$ such that almost-sure convergence fails to hold, even though $\E(|X_{\bm{k}}|^{r'})<\infty$ for all $r'\in (0,r)$. To see this, consider $X_{\i}=V_{i_1}$, where the $(V_i)_{i\geq 1}$ are i.i.d. $\alpha$-stable variables, with stability parameter $r \in (0,2)$ (and mean 0 when $r>1$). Then $E[|X_{\bm{k}}|^{r'}]<\infty$ for all $r'\in (0,r)$. On the other hand, considering here the case $r>1$ (the case $r\leq 1$ can be treated similarly),
$$\frac{1}{n^{k-1+1/r}}\sum_{\i\in\I_{n,k}}(X_{\i}-\E(X_{\bm{k}}))= \frac{(n-1)!}{(n-k)! n^{k-1}} \frac{1}{n^{1/r}}\sum_{i=1}^n V_i.$$
Then, because $n^{-1/r}\sum_{i=1}^n V_i\stackrel{d}{=}V_i$, $n^{-(k-1+1/r)}\sum_{\i\in\I_{n,k}}X_{\i}-\E[X_{\bm{k}}]$ does not converge to 0. Still, weaker moment conditions involving projections were shown to be sufficient for U-statistics \citep{teicher1998marcinkiewicz}. Similar conditions involving expectations conditional on unobserved terms $(U_{\{i\}})_{i\geq 1}$ are likely to be sufficient here as well.
\end{rem}

Theorem \ref{thm:MZL}, or rather a slight extension of it, implies the following law of iterated logarithm.

\begin{thm}
  Let $X$ be a dissociated, jointly exchangeable array such that $\E(X_{\bm{k}}^2)<\infty$ and let
  $$V:=\frac{1}{[(k-1)!]^2}\sum_{\substack{\pi\in\mathfrak{S}(\{1,...,k\})\\\pi'\in\mathfrak{S}(\{1,k+1,...,2k-1\})}} \Cov\left(X_{\pi(\bm{k})}, X_{\pi(1,k+1,...,2k-1)}\right).$$
  Then, $V\geq 0$ and
$$\lim\sup_{n\to\infty} \pm \frac{\sum_{\i\in \I_{n,k} }\left(X_{\i}-\E(X_{\bm{k}})\right)}{\sqrt{2 n^{2k-1}\log\log n}} = \sqrt{V} \; \text{a.s.}$$
\label{thm:LIL}
\end{thm}

The result improves upon that of \cite{scott1985law} by only requiring $\E(X_{\bm{k}}^2)<\infty$, exactly as with i.i.d. variables. Its proof can be summarized as follows (here again with $k=3$). First, we use Decomposition \eqref{eq:decomp2}. The term $\sum_{\i\in \I_{n,k}} H_1(X)_{\i}$ corresponds to a sum over i.i.d. terms, on which we apply the usual law of iterated logarithm. Next, it turns out that when $k\geq 2$, the proof of Theorem \ref{thm:MZL} extends to $r=2$ for the terms $\sum_{\i\in \I_{n,k}} H_2(X)_{\i}$ and $\sum_{\i\in \I_{n,k}} H_3(X)_{\i}$. As a result, these (properly normalized) terms are asymptotically negligible.

\section{Proof of Theorem \ref{thm:MZL} for $r>1$}
\label{sec:proof}

We recall that the proofs of Theorem \ref{thm:MZL} for $r<1$ and Theorem \ref{thm:LIL} are deferred to the supplementary material. We focus hereafter on $k\ge 2$ as the case $k=1$ is the usual result for i.i.d. data. Let us first introduce additional notation. Let $\bm{0}=(0,...,0)\in\N^k$ and $\bm{1}=(1,...,1)\in\N^k$. For any $\i=(i_1,...,i_k)$  and $\j=(j_1,...,j_k)$ in $\N^k$, $\i \odot \j$ stands for $(i_1 j_1,...,i_k j_k)$. With a slight abuse of notation, $\{\i\}$ denotes for any $\i=(i_1,...,i_k)\in \N^k$ the set of distinct elements of $(i_1,...i_k)$ and $\{\i\}^+ = \{\i\} \cap \N^+$. For any $\sigma\in\mathfrak{S}(\{1,...,k\})$ and $\i\in\I_k$, we let $\i_\sigma=(i_{\sigma(1)},...,i_{\sigma(k)})$. We say that $X$ is symmetric if for all $\i\in\I_k$ and $\sigma\in\mathfrak{S}(\{1,...,k\})$, $X_{\i_\sigma}=X_{\i}$.

\medskip
Next, for any $j\in\{0,...,k\}$, we let
$$\mathcal{E}_j=\left\{(e_1,...,e_k) \in \{0,1\}^k: \sum_{\ell=1}^k e_\ell = j\right\}, \quad \mathcal{E}=\cup_{j=0}^k\mathcal{E}_j.$$
For any $(\e,\e')\in\mathcal{E}^2$, we let $\e\leq \e'$ (resp. $\e\geq \e'$) when $\e_\ell\leq \e_\ell'$ (resp. $\e_\ell\geq \e_\ell'$) for $\ell\in\{1,...,k\}$. For any $\i\in\I_{n,j}$ and $\e\in\mathcal{E}_j$, we let $\i^{\e}$ be the $k$-dimensional vector with $i_1$ at the first non-null entry of $\e$, $i_2$ at the second non-null entry of $\e$ etc. and 0 elsewhere. For instance, if $k=4$, $j=2$, $\e=(0,1,0,1)$ and $\i=(5,3)$, we have $\i^{\e}=(0,5,0,3)$. Finally, for any $q\in\N^+$ and $A\subset \N^q$, we let $\overrightarrow{A}=\{\i \in \overline{A}: i_1< ... < i_q\}$. We also use $d$ as a shortcut for $k-1+1/r$.

In the following, we assume without loss of generality that $\E(X_{\bm{k}})=0$.

\subsection{Preliminary results}

The convergence in $L^r$ for $r>1$ relies on the following result, proved in \cite{davezies2021empirical}.

\begin{lem}[Symmetrization lemma]
	Let $r>1$ and $X$ be an exchangeable and dissociated array such that $\E\left(\left|X_{\bm{k}}\right|^r\right)<\infty$ and $\E\left(X_{\bm{k}}\right)=0$. Let $(\varepsilon_{A})_{A\subset \mathbb{N}}$ a family of independent Rademacher variables independent of $X$.  There exists $D_{r,k}$ a constant depending only on $r$ and $k$ and $k$ jointly exchangeable arrays $(X^j_{\i})_{\i\in \mathbb{I}_k}$ (for $j\in \{1,...,k\}$) such that $X^j_{\bm{k}}\stackrel{d}{=}X_{\bm{k}}$ and
	\begin{align*}
	\E\left[\left|\sum_{\i \in \mathbb{I}_{n,k}}X_{\i}\right|^r\right]&\leq D_{r,k}\sum_{j=1}^k\sum_{\e \in \mathcal{E}_{j}}\E\left(\left|\sum_{\i \in \mathbb{I}_{n,k}}\varepsilon_{\{\i\odot \e\}^+} X^j_{\i}\right|^r\right).
	\end{align*}
		\label{lem:sym}
\end{lem}

Almost-sure convergence relies on the AHK representation for the dissociated and jointly exchangeable array $(X_i)_{i\in\I_k}$. It states that there exists a function $\tau$ and i.i.d. uniform $[0,1]$ variables $(U_A)_{A\subset \cup_{r=1}^k \I_k}$ (and $U_{\emptyset}=1$ a.s.) such that $X_{\i}=\tau\left(\left(U_{\{\e\odot\i\}^+}\right)_{\e \in \mathcal{E}}\right)$. In the definition of $\tau$, the arguments are implicitly ordered by choosing a specific ordering on $\mathcal{E}$. We can now derive decompositions similar to the Hoeffding decomposition of U-statistics. Specifically, consider $Z_{\i}=f_{\i}(X_{\i})$ for any measurable maps $(f_{\i})_{\i\in\I_k}$. As soon as $\E|Z_{\i}|<\infty$ for all $\i\in\I_k$, we define, for all $\e\in \mathcal{E}$, $P_{\e}(Z)_{\i}=\E\left(Z_{\i}|(U_{\{\i\odot \e'\}^+})_{\e'\leq \e}\right)$. We then let  $H_0(Z)_{\i}=Q_{\zero}(Z)_{\i}=P_{\zero}(Z)_{\i}=\E(Z_{\i})_{\i}$ for any $\i\in \I_{n,k}$. For $\ell=1,...,k$, let
$$Q_{\e}(Z)_{\i}= \sum_{\e'\leq \e}(-1)^{\sum_{j=1}^ke_j-\sum_{j=1}^ke'_j} P_{\e'}(Z)_{\i}$$
and $H_{\ell}(Z)_{\i}= \sum_{\e \in \mathcal{E}_{\ell}} Q_{\e}(Z)_{\i}$. The following lemma gives three key properties on  $P_{\e}(Z)_{\i}, Q_{\e}(Z)_{\i}$ and $ H_{\ell}(Z)_{\i}$. In particular, Equality \eqref{eq:decompo}  may be seen as a  Hoeffding decomposition of $Z_{\i}$, as the components can be shown to be orthogonal in $L^2$. The proof of this lemma is given in the supplementary material.

\begin{lem}[Three properties of $P_{\e}(Z)_{\i}, Q_{\e}(Z)_{\i}$ and $H_{\ell}(Z)_{\i}$]\label{lem:prop_proj}
We have
	\begin{align}
	&Z_{\i}=\sum_{\ell=0}^{k}H_{\ell}(Z)_{\i}
	\label{eq:decompo},\\
	&H_{\ell}= \sum_{j=0}^{\ell}(-1)^{\ell-j}\binom{k-j}{\ell-j}\sum_{\e\in \mathcal{E}_j}P_{\e} \label{eq:H_P}.
	\end{align}
Moreover, for any symmetric array $X$, $\e\in \mathcal{E}_{\ell}$ (with $\ell\in\{0,...,k\}$), $\sigma\in \mathfrak{S}\{\bm{k}\}$ such that $\e_{\sigma}=\e$ and $\i\in \I_k$, we have $Q_{\e}(X)_{\i_\sigma}=Q_{\e}(X)_{\i}$.
\end{lem}


\subsection{Convergence in $L^r$}

Using Lemma~\ref{lem:sym} and the Khintchine inequality conditional on the $(X^j_{\i})_{\i\in \I_{k}}$, we obtain
\begin{align*}
\E\left[\left|\sum_{\i \in \mathbb{I}_{n,k}}X_{\i} \right|^r\right]&\leq D_{r,k}\sum_{j=1}^k\sum_{\e \in \mathcal{E}_{j}}\E\left(\left|\sum_{\i \in \mathbb{I}_{n,k}}\varepsilon_{\{\i\odot \e\}^+} X^j_{\i}\right|^r\right)\\
&= D_{r,k}\sum_{j=1}^k\sum_{\e \in \mathcal{E}_{j}}\E\left(\left|\sum_{\i \in \mathbb{I}_{n,j}}\varepsilon_{\{\i\}}\sum_{\i'\in \overline{(\{1,...,n\}\backslash\{\i\})^{k-j}}}X^j_{\i^{\e}+\i'^{\bm{1}-\e}} \right|^r\right)\\
&\leq D_{r,k} B_r \sum_{j=1}^k\sum_{\e\in \mathcal{E}_{j}}\E\left(\left|\sum_{\i \in \mathbb{I}_{n,j}}\left(\sum_{\i'\in \overline{(\{1,...,n\}\backslash\{\i\})^{k-j}}}X^j_{\i^{\e}+\i'^{\bm{1}-\e}}\right)^2\right|^{r/2}\right),
\end{align*}
where we recall that $\i^{\e}$ is defined at the beginning of the section and $B_r$ is a constant depending only on $r$. For each $j=1,...,k$ and $\e\in \mathcal{E}_j$, use $(a+b)^2\leq 2a^2+2b^2$ and $|a+b|^{r/2}\leq |a|^{r/2}+|b|^{r/2}$ to deduce:
\begin{align*}
&\E\left(\left|\sum_{\i \in \mathbb{I}_{n,j}}\left(\sum_{\i'\in \overline{(\{1,...,n\}\backslash\{\i\})^{k-j}}}X^j_{\i^{\e}+\i'^{\bm{1}-\e}}\right)^2\right|^{r/2}\right)\\
&\leq 2^{r/2}\E\left(\left|\sum_{\i \in \mathbb{I}_{n,j}}\left(\sum_{\i'\in \overline{(\{1,...,n\}\backslash\{\i\})^{k-j}}}X^j_{\i^{\e}+\i'^{\bm{1}-\e}}\mathds{1}\{|X^j_{\i^{\e}+\i'^{\bm{1}-\e}}|\leq M\}\right)^2\right|^{r/2}\right)\\
&+2^{r/2}\E\left(\left|\sum_{\i \in \mathbb{I}_{n,j}}\left(\sum_{\i'\in \overline{(\{1,...,n\}\backslash\{\i\})^{k-j}}}X^j_{\i^{\e}+\i'^{\bm{1}-\e}}\mathds{1}\{|X^j_{\i^{\e}+\i'^{\bm{1}-\e}}|> M\}\right)^2\right|^{r/2}\right).
\end{align*}
By Jensen's inequality $\E(|V|^{r/2})\leq \E(|V|)^{r/2}$, the first term can be bounded by $2^{r/2} n^{(k-j/2)r}M^{r}$. Using $|a+b|^{r/2}\leq |a|^{r/2}+|b|^{r/2}$, Jensen's inequality $\left(\sum_{a\in A} |V_a|\right)^r\leq \card(A)^{-1}\sum_{a\in A}\card(A)^r|V_a|^r$ and $\E\left(|X^j_{\i^{\e}+\i'^{\bm{1}-\e}}|^r\mathds{1}\{|X^j_{\i^{\e}+\i'^{\bm{1}-\e}}|> M\}\right)= \E\left(|X_{\bm{k}}|^r\mathds{1}\{|X_{\bm{k}}|> M\}\right)$, the second term can be bounded by $2^{r/2}n^{j}n^{(k-j)r}\E(|X_{\bm{k}}|^r\mathds{1}\{|X_{\bm{k}}|> M\})$. As a result,
\begin{align*}
& \frac{1}{n^{(k-1)r+1}}\E\left[\left|\sum_{\i \in \mathbb{I}_{n,k}}X_{\i} \right|^r\right] \\
&\leq D_{r,k}B_r2^{r/2}\sum_{j=1}^k \binom{k}{j} \left(n^{r(1-j/2)-1}M^r+n^{(j-1)(1-r)}\E\left(|X_{\bm{k}}|^r\mathds{1}\{|X_{\bm{k}}|>M\}\right)\right)\\
&\leq D_{r,k}B_r2^{r/2}(2^k -1)\left(n^{r/2-1}M^r+\E\left(|X_{\bm{k}}|^r\mathds{1}\{|X_{\bm{k}}|>M\}\right)\right).
\end{align*}
Considering $M$ sufficiently large to ensure that $\E\left(|X_{\bm{k}}|^r\mathds{1}\{|X_{\bm{k}}|> M\}\right)$ is arbitrarily small, we deduce that $\lim\sup_{n}\E\left[\left|n^{-d}\sum_{\i \in \I_{n,k}}X_{\i} \right|^r\right]=0$.

\subsection{Almost-sure convergence}

Up to considering $Z_{\i}=\sum_{\sigma\in\mathfrak{S}(\{\bm{k}\})} [X_{\i_\sigma} - \E(X_{\bm{k}})]/k!$ instead of $X_{\i}$, we can assume without loss of generality that $\E(X_{\bm{k}})=0$ and that $X$ is symmetric. Then:
$$\sum_{\i \in \I_{n,k}}X_{\i}=k!\sum_{\i \in \overrightarrow{\I_{n,k}}}X_{\i}=k!\sum_{\ell=1}^{k}\sum_{\i \in \overrightarrow{\I_{n,k}}}H_{\ell}(X)_{\i},$$
where we recall the notation $\overrightarrow{A}=\{\i\in\overline{A}: i_1<...<i_q\}$. Following the U-statistic's terminology, we say that $X$ is degenerate of order $\ell-1$ if $H_{\ell}(X)_{\i}$ is not constant whereas for all $\ell'\in \{1,...,\ell-1\}$, $H_{\ell'}(X)_{\i}$ is constant (and then equal to 0 as we  assumed $\E(X_{1,...k})=0$).

\paragraph{Degenerate arrays of order $k-1$.}
We first assume that $X$ is degenerate of order $k-1$, that is $X_{\i}=H_{k}(X)_{\i}$. By Kronecker’s lemma, it is sufficient to show that
$$\sum_{j=k}^{n}\frac{1}{j^{d}}\sum_{\i \in \overrightarrow{\I_{j-1,k-1}}} X_{\i,j} \to 0 \quad \text{a.s.}$$
Let $T_j=j^{-d}\sum_{\i \in \overrightarrow{\I_{j-1,k-1}}} X_{\i,j}$ and $\mathcal{F}_j=\sigma\left(U_{A};A\subset\{1,...,j\},\card(A)\leq k\right)$. Since $X$ is degenerate of order $k-1$, we have $\E\left(X_{\bm{k}}|\left(U_{A}\right)_{A\subset\{1,...,k-1\}}\right)=0$. Then, $\E\left(T_j|\mathcal{F}_{j-1}\right)=0$, meaning that $(T_j,\mathcal{F}_j)$ is a martingale difference. Hence, $\left(\sum_{j=\ell}^nT_j,\mathcal{F}_n\right)$ is a martingale. We then just have to prove that $\sup_{n}\E\left|\sum_{j=k}^nT_j\right|<\infty$ \citep[see for instance Chapter 10, Theorem 12.2 in][] {gut2013probability}. To this end, we use
$\E\left|\sum_{j=k}^nT_j\right|\leq A_n + B_n$, with
$A_n =  \E\left|\sum_{j=k}^nj^{-d}\sum_{\i\in\overrightarrow{\I_{j-1,k-1}}}X_{\i,j}\right.$ $\left.\mathds{1}\{|X_{\i,j}|\leq j^{d}\}\right|$ and $B_n = \E\left|\sum_{j=k}^nj^{-d}\sum_{\i\in\overrightarrow{\I_{j-1,k-1}}}X_{\i,j} \mathds{1}\{|X_{\i,j}|> j^{d}\}\right|$.
We show that the suprema of both terms are finite. First,
\begin{align}
B_n &\leq\sum_{j=k}^nj^{-d}\frac{(j-1)!}{(j-k)!}\E\left(\left|X_{\bm{k}}\right|\mathds{1}\{|X_{\bm{k}}|> j^{d}\}\right)\label{eq:ineg1} \\
&\leq\E\left(\left|X_{\bm{k}}\right|\sum_{j=k}^nj^{-1/r}\mathds{1}\{|X_{\bm{k}}|> j^{d}\}\right)\notag \\
&\leq\E\left(\left|X_{\bm{k}}\right|\mathds{1}\{|X_{\bm{k}}|> (k-1)^{d}\}\int_{k-1}^{|X_{\bm{k}}|^{1/d}}t^{-1/r}dt\right)\notag \\
&\leq\E\left(\left|X_{\bm{k}}\right|\mathds{1}\{|X_{\bm{k}}|> (k-1)^{d}\}\left[\frac{t^{1-1/r}}{1-1/r}\right]_{k-1}^{|X_{\bm{k}}|^{1/d}}\right)\notag \\
&\leq \frac{r}{r-1}\E\left(\left|X_{\bm{k}}\right|^{1+\frac{1-1/r}{d}}\right) \notag \\
&\leq \frac{r}{r-1}\E\left(\left|X_{\bm{k}}\right|^{r}\right)  \label{eq:ineg2},
\end{align}
where the last inequality follows using $1+(1-1/r)/d\leq r$. Hence, $\sup_n B_n <\infty$. 

\medskip
Now let us turn to $A_n$. First, applying Decomposition \eqref{eq:decompo} to $X_{\i,j}\mathds{1}\{|X_{\i,j}|\leq j^{d}\}$ for a fixed $j$ and using the fact that $H_{\ell}(X)=0$ for $\ell=0,...,k-1$, we get, for any $\i\in \I_{j-1,k-1}$,
\begin{align*}
X_{\i,j}\mathds{1}\{|X_{\i,j}|\leq j^{d}\}
&=H_k(X\mathds{1}\{|X|\leq j^{d}\})_{\i,j}+H_0(X\mathds{1}\{X\leq j^{d}\})_{\i,j}+\sum_{\ell=1}^{k-1}H_{\ell}(X\mathds{1}\{|X|\leq j^{d}\})_{\i,j}\\
&=H_k(X\mathds{1}\{|X|\leq j^{d}\})_{\i,j}+H_0(X\mathds{1}\{X> j^{d}\})_{\i,j}+\sum_{\ell=1}^{k-1}H_{\ell}(X\mathds{1}\{|X|> j^{d}\})_{\i,j}.
\end{align*}
Let $A_n=A_{1n}+ A_{2n}+A_{3n}$ be the decomposition of $A_n$ corresponding to the first, second and third terms in the previous display. We have
\begin{align}
A_{2n}&\leq \sum_{j=k}^nj^{-d}\frac{(j-1)!}{(j-k)!}\E\left(\left|X_{\bm{k}}\right|\mathds{1}\{|X_{\bm{k}}|> j^{d}\}\right)\notag\\
&\leq \frac{r}{r-1}\E\left(\left|X_{\bm{k}}\right|^r\right)\label{eq:A2},
\end{align}
where \eqref{eq:A2} uses the inequalities from \eqref{eq:ineg1} to \eqref{eq:ineg2}. Thus, $\sup_n A_{2n}<\infty$. Let us turn to $A_{3n}$. First, for any $\e\in \mathcal{E}$, we have
\begin{align}
&\E\left|\sum_{j=k}^{n}j^{-d}\sum_{\i\in \overrightarrow{\I_{j-1,k-1}}}P_{\e}(X\mathds{1}\{X> j^{d}\})_{\i,j}\right|\notag\\
&\leq \sum_{j=k}^{n}j^{-d}\sum_{\i\in \overrightarrow{\I_{j-1,k-1}}}\E\left(\E\left(|X_{\i,j}|\mathds{1}\{X> j^{d}\})_{\i,j}|U_{\{(\i,j)\odot \e\}^+}\right)\right)\notag\\
&\leq \sum_{j=k}^{n}j^{-d}\frac{(j-1)!}{(j-k)!}\E\left(\left|X_{\bm{k}}\right|\mathds{1}\{\left|X_{\bm{k}}\right|>j^{d}\}\right)\notag\\
&\leq \frac{r}{r-1}\E\left(\left|X_{\bm{k}}\right|^r\right).\label{eq:A3}
\end{align}
Then, in view of \eqref{eq:H_P}, we have $\sup_n A_{3n}<\infty$. Finally, let us consider $A_{1n}$. We can write
\begin{align}
A_{1n}^2 & = \E\left[\sum_{j=k}^{n}j^{-d}\sum_{\i\in \overrightarrow{\I_{j-1,k-1}}}H_k(X\mathds{1}\{|X|\leq j^{d}\})_{\i,j} \right]^2\notag\\
& \leq \E\left[\left(\sum_{j=k}^{n}j^{-d}\sum_{\i\in \overrightarrow{\I_{j-1,k-1}}}H_k(X\mathds{1}\{|X|\leq j^{d}\})_{\i,j}\right)^2\right]\notag\\
&=\sum_{j=k}^{n}j^{-2d}\sum_{\i\in \overrightarrow{\I_{j-1,k-1}}}\E\left(H_k(X\mathds{1}\{|X|\leq j^{d}\})^2_{\i,j}\right)\notag \\
&\leq \sum_{j=k}^{n} j^{-2d+k-1}\E\left(X_{\bm{k}}^2\mathds{1}\{|X_{\bm{k}}|\leq j^{d}\}\right)\notag\\
&\leq \E\left[X_{\bm{k}}^2\left( \ind{k> |X_{\bm{k}}|^{1/d}} \int_{k-1}^{\infty}t^{-k+1-2/r}dt +\ind{k\le |X_{\bm{k}}|^{1/d}} \int_{|X_{\bm{k}}|^{1/d}-1}^{\infty}t^{-k+1-2/r}dt \right) \right] \notag\\
&\leq \frac{1}{k-2+2/r}\left\{k^{2d} + \E\left[X_{\bm{k}}^2 \left(|X_{\bm{k}}|^{1/d}-1\right)^{-k+2-2/r}\ind{|X_{\bm{k}}|^{1/d}\ge k}\right]\right\},\label{eq:A1}
\end{align}
where the second equality follows since the variables $(H_k(X\mathds{1}\{|X|\leq j^{d}\})_{\i,j})$ are uncorrelated. Now, note that for all $x\ge 0$,
$$x^2 \left(x^{1/d}-1\right)^{-k+2-2/r}\ind{x^{1/d}\ge k} \le \left(1-1/k\right)^{k-2d} x^{k/d}\le \left(1-1/k\right)^{k-2d} x^r.$$
Therefore, 
$$\E\left[X_{\bm{k}}^2 \left(|X_{\bm{k}}|^{1/d}-1\right)^{-k+2-2/r}\ind{|X_{\bm{k}}|^{1/d}\ge k}\right]\leq \left(1-1/k\right)^{k-2d} \E\left[X_{\bm{k}}^r\right]<\infty.$$
Hence, $\sup_n A_{1n}<\infty$, which completes the proof when $X$ is  degenerate of order $k-1$.

\paragraph{Other cases.}

In view of \eqref{eq:decompo}, it suffices to prove that for all $\ell\in\{0,...,k\}$, we have
\begin{equation}
	\label{eq:cv_as_H}
	\frac{1}{n^{d}}\sum_{\i \in \overrightarrow{\I_{n,k}}}H_{\ell}(X)_{\i}\to 0 \; \text{a.s.}
\end{equation}
First, note that $\sum_{\i \in \overrightarrow{\I_{n,k}}}H_{\ell}(X)_{\i}=\frac{1}{k!}\sum_{\i \in \I_{n,k}}H_{\ell}(X)_{\i}=\sum_{\e \in \mathcal{E}_{\ell}}\frac{1}{k!}\sum_{\i \in \I_{n,k}}Q_{\e}(X)_{\i}.$
By Lemma \ref{lem:prop_proj}, $Q_{\e}(X)_{\i_\sigma}=Q_{\e}(X)_{\i}$ for all $\sigma$ such that $\e_\sigma=\e$. Moreover, for each $\e\in \mathcal{E}_{\ell}$, $Q_{\e}(X)_{\i}$ only depends on $\i\odot\e$ rather than all components of $\i$. It follows that
$\sum_{\i \in \I_{n,k}}Q_{\e}(X)_{\i}=\frac{(n-\ell)!}{(n-k)!}\sum_{\i \in \I_{n,\ell}}R^{\e}_{\i},$
where $(R^{\e}_{\i})_{\i\in \I_{n,\ell}}$ is a symmetric $\ell$-dimensional jointly exchangeable array degenerate of order $\ell-1$ such that $\E\left(\left|R^{\e}_{1,...,\ell}\right|^r\right)<\infty$. It follows from the  previous paragraph that
$$\frac{1}{n^{\ell-1+1/r}}\sum_{\i \in \overrightarrow{\I_{n,\ell}}}R^{\e}_{\i}=\frac{1}{n^{\ell-1+1/r} \times \ell!}\sum_{\i \in \I_{n,\ell}}R^{\e}_{\i} \to 0 \;\text{a.s.}$$
Because $(n-\ell)!/(n-k)!\sim n^{k-\ell}$, we finally obtain \eqref{eq:cv_as_H}.

\bibliography{../biblio}

\newpage
\setcounter{page}{1}
\appendix
\begin{center}{\Huge Supplementary material}\end{center}

\medskip
\section{Proof of Lemma \ref{lem:prop_proj}}

First,
\begin{align*}
\sum_{\ell=0}^{k}H_{\ell}(Z)_{\i}&=\sum_{\e \in \mathcal{E}}\sum_{\e'\leq \e}(-1)^{\sum_{j=1}^ke_j-\sum_{j=1}^ke'_j}P_{\e'}(Z)_{\i}\\
&=\sum_{\e' \in \mathcal{E}}\sum_{\e'\leq \e}(-1)^{\sum_{j=1}^ke_j-\sum_{j=1}^ke'_j}P_{\e'}(Z)_{\i}\\
&=\sum_{\e' \in \mathcal{E}}\left(\sum_{j=0}^{k-\sum_{j'}e'_{j'}}(-1)^{j}\binom{k-\sum_{j'}e_{j'}}{j}\right)P_{\e'}(Z)_{\i}\\
&=\sum_{\e' \in \mathcal{E}}(1-1)^{k-\sum_{j'}e'_{j'}}P_{\e'}(Z)_{\i}\\
&=P_{\un}(Z)_{\i}=Z_{\i}.
\end{align*}
Next,
\begin{align*}
H_{\ell}&=\sum_{\e\in \mathcal{E}_{\ell}}\sum_{\e'\leq \e}(-1)^{\ell-\sum_{m}e'_m} P_{\e'}\\
&=\sum_{j=0}^{\ell}\sum_{\e'\in \mathcal{E}_j}(-1)^{\ell-j}\binom{k-j}{\ell-j}P_{\e'}\\
&=\sum_{j=0}^{\ell}(-1)^{\ell-j}\binom{k-j}{\ell-j}\sum_{\e\in \mathcal{E}_j}P_{\e}.
\end{align*}
To prove the last claim, note that if $X$ is symmetric, we have for all $\i\in\I_k$,
\begin{align*}
Q_{\e}(X)_{\i_{\sigma}}&=\sum_{\e'\leq \e}\left(-1\right)^{\ell-\sum_{m=1}^ke'_m}\E\left(X_{\i_{\sigma}}|(U_{\i_{\sigma} \odot \e''})_{\e''\leq \e'}\right)\\
&=\sum_{\e'\leq \e}\left(-1\right)^{\ell-\sum_{m=1}^ke'_m}\E\left(X_{\i}|(U_{\i_{\sigma} \odot \e''})_ {\e''\leq \e'}\right)\\
&=\sum_{\e'\leq \e}\left(-1\right)^{\ell-\sum_{m=1}^ke'_m}\E\left(X_{\i}|(U_{\i \odot \e''_{\sigma^{-1}}})_ {\e''\leq \e'}\right)\\
&=\sum_{\e'\leq \e}\left(-1\right)^{\ell-\sum_{m=1}^ke'_m}\E\left(X_{\i}|(U_{\i \odot \e''_{\sigma^{-1}}})_ {\e''\leq \e'_{\sigma}}\right)\\
&=\sum_{\e'\leq \e}\left(-1\right)^{\ell-\sum_{m=1}^ke'_m}\E\left(X_{\i}|(U_{\i \odot \e''_{\sigma^{-1}}})_ {\e''_{\sigma^{-1}}\leq \e'}\right)\\
&=\sum_{\e'\leq \e}\left(-1\right)^{\ell-\sum_{m=1}^ke'_m}\E\left(X_{\i}|(U_{\i \odot \e''})_{\e''\leq \e'}\right)\\
&=Q_{\e}(X)_{\i}.
\end{align*}

\section{Proof of Theorem \ref{thm:MZL} for $r<1$}

\subsection{Convergence in $L^r$}

An inspection of the proof of Theorem 10.3 in \citep[][p. 311]{gut2013probability} shows that it does not rely on the independence of the variables but only on their identical distribution. Hence, it directly applies here.

\subsection{Almost-sure convergence}

An inspection of the proof of Theorem 1 in \cite{Gine1992} shows that it only exploits the fact that the individual variables are identically distributed. Hence, it also applies here.

\section{Proof of Theorem \ref{thm:LIL}}

Again, we focus on $k\ge 2$ as the case $k=1$ is the usual result for i.i.d. data. We also assume without loss of generality that $\E\left(X_{\bm{k}}\right)=0$. Let us first prove the result when $X$ is symmetric. The inequalities \eqref{eq:ineg2}, \eqref{eq:A2} and \eqref{eq:A3} actually hold for $r=2$. Also, Inequality \eqref{eq:A1} is valid for $r=2$. It follows that $n^{-(k-1/2)}\sum_{\i \in \I_{n,k}}H_{\ell}(X)_{\i}$ converges almost surely to 0 for $\ell\geq 2$. Then, by Decomposition \eqref{eq:decompo},
\begin{align*}
\frac{1}{n^{k-1/2}}\sum_{\i \in \I_{n,k}}X_{\i}&=\frac{1}{n^{k-1/2}}\sum_{\i \in \I_{n,k}}H_{1}(X)_{\i}+o_{\text{a.s.}}(1) \\
&=\frac{1}{n^{k-1/2}}\frac{(n-1)!}{(n-k)!}k\sum_{i=1}^n\E\left(X_{(i,n+1,...,n+k-1)}|U_{\{i\}}\right)\\
&=\left(\frac{k}{n^{1/2}}\sum_{i=1}^n\E\left(X_{(i,n+1,...,n+k-1)}|U_{\{i\}}\right)\right)\left(1+o(1)\right).\end{align*}
The usual law of iterated logarithm ensures that
\begin{align}\limsup_{n}\pm\frac{\sum_{\i \in \I_{n,k}}X_{\i}}{\sqrt{2n^{2k-1}\log\log(n)}}=k\sqrt{\V(\E\left(X_{\bm{k}}|U_{\{1\}}\right)}.\label{eq:LIL1}\end{align}
From the AHK decomposition, we have $\Cov\left(X_{\bm{k}}, X_{(1,k+1,...,2k-1)}|U_{\{1\}}\right)=0$. Next, symmetry of $X$ and the law of total covariance yield:
\begin{align}
V & = k^2\Cov\left(X_{\bm{k}}, X_{(1,k+1,...,2k-1)}\right)\notag\\
&=k^2\Cov\left(\E\left(X_{\bm{k}}|U_{\{1\}}\right), \E\left(X_{(1,k+1,...,2k-1)}|U_{\{1\}}\right)\right)\notag\\
&=k^2\V\left(\E\left(X_{\bm{k}}|U_{\{1\}}\right)\right) (\geq 0).\label{eq:LIL2}
\end{align}
When $X$ is not symmetric, we just replace $X_{\i}$ by $\sum_{\pi \in \mathfrak{S}(\{\i\})}X_{\pi(\i)}/k!$ in \eqref{eq:LIL1} and \eqref{eq:LIL2}.

\end{document}